\documentclass[12pt]{article}
\usepackage{amsmath}
\usepackage{amssymb}
\usepackage{amscd}
\usepackage{epsfig}

\def\no{\noindent}
\numberwithin{equation}{subsection}
\numberwithin{figure}{subsection}

\newenvironment{theorem}{\par\bigskip\noindent\refstepcounter{equation}{\bf Theorem \theequation.}\em}{\em\par\bigskip\noindent}

\newenvironment{proposition}{\par\bigskip\noindent\refstepcounter{equation}{\bf Proposition \theequation.}\em}{\em\par\bigskip\noindent}

\newenvironment{conjecture}{\par\bigskip\noindent\refstepcounter{equation}{\bf Conjecture \theequation.}\em}{\em\par\bigskip\noindent}
\newenvironment{definition}{\par\bigskip\noindent\refstepcounter{equation}{\bf Definition \theequation.}}{\par\bigskip\noindent}

\newenvironment{question}{\par\bigskip\noindent\refstepcounter{equation}{\bf Question \theequation.}}{\par\bigskip\noindent}

\newcommand{\N}{\mathbb N}
\newcommand{\R}{\mathbb R}
\newcommand{\T}{\mathbb T}

\newcommand{\Z}{\mathbb Z}

\newcommand{\qed}{{\unskip\nobreak\hfil
        \penalty50\hskip1em\hbox{}\nobreak\hfil
        $\square$\parfillskip=0pt\finalhyphendemerits=0 \par}}
\newcommand{\proof}{\no{\em Proof.\ }}

\def\sys {{\it sys}} \def\stsys {{\it stsys}} \def\confsys {{\it
confsys}} \def\pisys {{\it sys}\pi} \def\phisys {\phi{\it sys}}
\def\vol{{\it vol}} \def\area{{\it area}} \def\morphism{{\it
morphism}} \def\diam{{\it diam}}

\def\length{{\it length}}
\def\len{{\it len}}

\def\C {{\mathbb C}} 

\def\CC{{the first author\ }} 
\def\MK{{the second author\ }} 

\def\ie {{\it i.e.\ }} 
\def\eg {{\it e.g.\ }} 
\def\cf {\hbox{\it cf.\ }}

\def\g {{\bf {g}}} 
\def\gen {{\it s}} 
\def\sing {{\sigma}}

\def\3P{{\it 3P}} \def\AJ{{\it AJ}} \def\XX{{\overline{X}}}

\def \SGM{{\it SGM}} 
\def \MLS{{\it MLS}}

\def \SR {{\it SR}}

\def\rat{{\angle\!^\mid \! \left(\frac{x}{2}, \frac{\pi}{8}\right)\ }}
\def\fillrad{{\it FillRad}}

\def\cstok#1{\leavevmode\thinspace\hbox{\vrule\vtop{\vbox{\hrule\kern1pt
    \hbox{\vphantom{\tt/}\thinspace{\tt#1}\thinspace}}
      \kern1pt\hrule}\vrule}\thinspace}

\def\stm{
\begin{figure}
\[
\begin{tabular}[t]{|
@{\hspace{3pt}}p{.9in}|| @{\hspace{3pt}}p{2.4in}|
@{\hspace{3pt}}p{2.3in}| } \hline & \cstok{1} $k=1$ & \cstok{2} $k\ge
2$ \\ \hline\hline \cstok{A}\break homotopy\break $k$-systole &
Loewner's inequality for $\T^2$ \eqref{(1.1)}, Pu's for $\R P^2$
\eqref{pu3}, Gromov's for essential $X$ \eqref{sys1} & \\ \hline
\cstok{B}\break homology\break $k$-systole & if cuplength$(X^n)=n$
then Gromov proves inequality \eqref{525} & freedom reigns: Gromov
\eqref{44}, Babenko \eqref{41}, \eqref{ib}; Katz, Suciu \eqref{41};
freedom of {$S^1\!\!\times\!  S^2$ over $\Z_2$}: Freedman \eqref{4.2}
\\ \hline \cstok{C}\break stable homology $k$-systole & sharp
inequality if $b_1(X^n)$ = cup\-length$(X^n)=n$ (Gromov, based on
Burago-Ivanov) \eqref{(1.4)} & multiplicative relations in $H^*(X)$
entail inequalities: Gromov, Hebda \eqref{hebda}, Bangert-Katz
\eqref{517}, \eqref{(5.2)} \\ \hline \cstok{D}\break interpolation
homotopy--homology & ``fiberwise'' inequality if Abel-Jacobi map
``surjective'' (Katz-Kreck-Suciu) \cite{KKS}, \cf \eqref{ptl}& \\
\hline \cstok{E}\break conformal $k$-systole & Accola, Blatter
(section \ref{highergenus}); for a surface $\Sigma_\gen$, logarithmic
in genus~$\gen$: Buser-Sarnak \eqref{bs} & $X^4$, indefinite case:
polynomial in $\chi(X)$, modulo surjectivity of period map: Katz
\eqref{T} \\ \hline \cstok{F}\break contractible $k$-systole & short
contractible geodesics: Croke \eqref{214}, Maeda, Nabutovsky and
Rotman \eqref{nrs}, Sabourau \eqref{65} & \\ \hline
\end{tabular}
\]
\caption{\textsf{A 2-D map of systolic geometry}}
\label{fig1}
\end{figure}
}

\begin{document}

\author{Christopher B. Croke\thanks{Supported by NSF grant DMS
02-02536.}\ ~and Mikhail G. Katz\thanks{Supported by ISF grant no.\
620/00-10.0}}

\date{\today} 

\title{Universal volume bounds in Riemannian manifolds\thanks{To
appear in Surveys in Differential Geometry (2003)}}

\maketitle
\begin{abstract}
In this survey article we will consider universal lower bounds on the
volume of a Riemannian manifold, given in terms of the volume of lower
dimensional objects (primarily the lengths of geodesics).  By
`universal' we mean without curvature assumptions.  The restriction to
results with no (or only minimal) curvature assumptions, although
somewhat arbitrary, allows the survey to be reasonably short.
Although, even in this limited case the authors have left out many
interesting results.
\end{abstract}

\tableofcontents

\setcounter{section}{1}
\setcounter{subsection}{0}

\subsection{Introduction}

In the present article, we will consider $n$-dimensional Riemannian
manifolds $(X^n,\g)$, for the most part compact.  We plan to survey
the results and open questions concerning volume estimates that do not
involve curvature (or involve it only very weakly).  Our emphasis will
be on open questions, and we intend the article to be accessible to
graduate students who are interested in exploring these questions.  In
choosing what to include here, the authors have concentrated on
results that have influenced their own work and on recent developments.  
In particular, we are
only able to mention some of the highlights of M.~Gromov's seminal
paper~\cite{Gr1}.  His survey devoted specifically to systoles
appeared in \cite{Gr2}.  The interested reader is encouraged to
explore those papers further, as well as his recent book \cite{Gr3}.

The estimates we are concerned with are lower bounds on the Riemannian
volume (\ie the $n$-dimensional volume), in terms of volume-minimizing
lower dimensional objects.  For example, one lower dimensional volume
we will consider will be the infimum of volumes of representatives of
a fixed homology or homotopy class.  This line of research was,
apparently, originally stimulated by a remark of Ren\'e Thom, in a
conversation with Marcel Berger in the library of Strasbourg
University in the 1960's, not long after the publication of
\cite{Ac60,Bl61}.  Having been told of the latter results (discussed
in sections \ref{highergenus} and \ref{51bis}), Professor Thom
reportedly exclaimed: ``Mais c'est fondamental! [These results are of
fundamental importance!].''  Ultimately this lead to the so called
isosystolic inequalities, \cf section \ref{katz}.  An intriguing
historical discussion appears in section ``Systolic reminiscences''
\cite[pp.~271-272]{Gr3}.

We will devote section \ref{2dim} to the study of lower bounds on the
area of surfaces (two dimensional manifolds) in terms of the length of
closed geodesics.  This is where the subject began, with the theorems
of Loewner and Pu, and where the most is known.  

In section \ref{gromov} we discuss some of the results and questions
that come from \cite{Gr1}. In particular, inequality
\eqref{filradvol}, which provides a lower bound for the total volume
in terms of an invariant called the Filling Radius, is one of the main
tools used to obtain isosystolic inequalities in higher dimensions.
We also discuss Gromov's notion of Filling Volume, as well as some of
the open questions and recent results relating the volume of a compact
Riemannian manifold with boundary, to the distances among its boundary
points \cite{B-C-G,Cr1, C-D-S,C-D, Iv1}.

In section \ref{katz}, we see a collection of results that show that
the inequalities envisioned in M. Berger's original question are often
violated.  Unless we are dealing with one-dimensional objects, such
isosystolic inequalities are systematically violated by suitable
families of metrics.  This line of work was stimulated by M.~Gromov's
pioneering example \eqref{44}, \cf \cite{Gr2, Ka1, Pit, BabK, BKS,
KS1, KS2, Bab, Ka2}.

Section \ref{stablesystole} discusses how the systolic inequalities
that we saw failed in Section \ref{katz}, can in fact hold if we pass
to stable versions.  We survey the known results of this nature
\cite{Gr1, He, BanK1, BanK2}, and also examine the related conformal
systolic invariants and their asymptotic behavior, studied in
\cite{BS, Ka3}.

In section \ref{inject}, we discuss bounds on the volume in terms of
the injectivity radius of a compact Riemannian manifold.  We discuss
the (sharp) isoembolic inequality \eqref{bergerisoemb} of M. Berger,
as well as the local version, \ie volumes of balls, by Berger and \CC.
We survey some of the extensions of these results, with an emphasis on
the open questions and conjectures.

\subsection{Area and length of closed geodesics in 2 dimensions}
\label{2dim}

In this section, we will restrict ourselves to 2-dimensional
Riemannian manifolds $(X^2,\g)$ and discuss lower bounds on the total
area, $area(\g)$, in term of the least length of a closed geodesic
$L(\g)$.  In any dimension, the shortest loop in every nontrivial
homotopy class is a closed geodesic.  We will denote by
$\pisys_1(\g)$, the least length (often referred to as the ``Systole''
of $\g$) of a noncontractible loop $\gamma$ in a compact,
non-simply-connected Riemannian manifold~$(X,\g)$:
\begin{equation}
\label{ps1}
\pisys_1(\g)=\min_{[\gamma]\not= 0\in \pi_1(X)}\length(\gamma).
\end{equation}
Hence $\pisys_1(\g)\geq L(\g)$.  In this section (with the exception
of section \ref{2-sphere} where there is no homotopy) we will consider
isosystolic inequalities of the following form: ${\pisys_1 (\g) ^2 \le
const\; \area(\g)}$.  This leads naturally to the notion of the
systolic ratio of an $n$-dimensional Riemannian manifold $(X,\g)$,
which is defined to be the scale-invariant quantity $\frac{\pisys_1
(\g)^n} {\vol_n(\g)}$.  We also define the {\em optimal systolic
ratio} of a manifold $X$ to be the quantity
\begin{equation}
\label{22}
\sup_\g\frac{\pisys_1(\g)^n}{\vol_n(\g)},
\end{equation}
where the supremum is taken over all Riemannian metrics $\g$ on the
given manifold $X$.  Note that in the literature, the reciprocal of
this quantity is sometimes used instead.  In this section, we will
discuss the optimal systolic ratio for $\T^2$, $\R P^2$ (section
\ref{loewnerpu}), and the Klein bottle $\R P^2\#\R P^2$ (section
\ref{optimal}).  These are the only surfaces for which the optimal
systolic ratio is known.  However, upper bounds for the systolic
ratio, namely isosystolic inequalities, are available for all compact
surfaces (section \ref{highergenus}) except of course for $S^2$ (but
see section \ref{2-sphere}).  Figure \ref{stmap} contains a chart
showing the known upper and lower bounds on the optimal systolic
ratios of surfaces.

\begin{figure}
\renewcommand{\arraystretch}{1.8}
\[
\begin{tabular}[t]{|
@{\hspace{3pt}}p{1.2in}|| @{\hspace{3pt}}p{1.9in}|
@{\hspace{3pt}}p{1.1in}| @{\hspace{3pt}}p{1.1in}| } \hline &
$\displaystyle\sup_\g\frac{\pisys_1(\g)^2}{\area(\g)}$ & numerical
value & where to find it \\ \hline\hline $X=\R P^2$ &
$\displaystyle =\frac{\pi}{2}$ \cite{Pu, Ta} & $\approx 1.5707$ & \eqref{pu3} \\
\hline infinite $\pi_1(X)$ & $\displaystyle <\frac{4}{3}$ \cite{Gr1}&
$<1.3333$ &\eqref{tqi} \\ \hline $X=\T^2$ &
$\displaystyle =\frac{2}{\sqrt{3}}$ (Loewner) & $\approx 1.1547$ & \eqref{(1.1)}
\\ \hline $X=\R P^2 \# \R P^2$ & $\displaystyle =\frac{\pi}{2^{3/2}}$
\cite{Bav86,Sa} & $\approx 1.1107$ &\eqref{2p} \\ \hline $X$ of genus 2 &
$>\left(27-18\sqrt{2}\right)^{-\frac{1}{2}}$ & $> 0.8047$ & \eqref{37}
\\ \hline $X$ of genus 3 & $\displaystyle \geq \frac 8 {7\sqrt{3}}$
\cite{Ca96} & $> 0.6598$ & section \ref{optimal} 
\\ \hline $X$ of genus 4 & $\displaystyle > \frac{9\sqrt 7}{40} $
\cite{Ca03,Ca3} & $> 0.5953$ & section \ref{optimal} 
\\ \hline $X$ of genus
$\gen$ & $\displaystyle< \frac {64} {4\sqrt{\gen} + 27}$ \cite{Gr1,
Ko} && \eqref{1.1} \\\hline
\end{tabular}
\]
\renewcommand{\arraystretch}{1}
\caption{Values for optimal systolic ratio of surface $(X,\g)$, in
decreasing order}
\label{stmap}
\end{figure}

\subsubsection{Inequalities of Loewner and Pu}
\label{loewnerpu}
The first results of this type were due to C. Loewner and P. Pu.
Around 1949, Carl Loewner proved the first systolic inequality, \cf
\cite{Pu}.  He showed that for every Riemannian metric $\g$ on the
torus $\T^2$, we have
\begin{equation}
\label{(1.1)}
\pisys_1(\g)^2\le\frac{2}{\sqrt{3}}\area(\g),
\end{equation} 
while a metric satisfying the boundary case of equality in
\eqref{(1.1)} is necessarily flat, and is homothetic to the quotient
of $\C$ by the lattice spanned by the cube roots of unity.

Two distinct optimal generalisations of \eqref{(1.1)} are available,
\cf \eqref{(5.2)} and \eqref{(1.4)}.  It follows from Gromov's
estimate \eqref{1.1} that aspherical surfaces satisfy Loewner's
inequality \eqref{(1.1)} if the genus is bigger than 50.  It is an
open question whether the genus assumption can be removed, but for
genus up to 30, some information can be deduced from the Buser-Sarnak
inequality \eqref{bs}.

We give a slightly modified version of M.~Gromov's proof \cite{Gr2},
using conformal representation and Cauchy-Schwartz, of Loewner's
theorem for the $2$-torus.  We present the following slight
generalisation: there exists a pair of closed geodesics on
$(\T^2,\g)$, of respective lengths $\lambda_1$ and $\lambda_2$, such
that
\begin{equation}
\label{23}
\lambda_1 \lambda_2 \leq \tfrac{2}{\sqrt{3}}\area(\g),
\end{equation}
and whose homotopy classes form a generating set for $\pi_1(\T^2)=
\Z\times\Z$.

\proof The proof relies on the conformal representation $\phi: \T_0
\rightarrow (\T^2,\g)$, where $\T_0$ is flat.  Here $\phi$ may be
chosen in such a way that $(\T^2,\g)$ and $\T_0$ have the same area.
Let $f$ be the conformal factor of $\phi$.

Let $\ell_0$ be any closed geodesic in $\T_0$.  Let $\{\ell_s\}$ be
the family of geodesics parallel to $\ell_0$.  Parametrize the family
$\{\ell_s\}$ by a circle $S^1$ of length $\sigma$, so that $\sigma
\ell_0 = \area(\T_0)$.  Thus $\T_0\to S^1$ is a Riemannian submersion.
Then $\area(\T^2)=\int_{\T_0} Jac_\phi = \int_{\T_0}f^2$.  By Fubini's
theorem, $\area(\T^2)=\int_{S^1}ds \int_{\ell_s} f^2 dt$.  By the
Cauchy-Schwartz inequality,
\[
\area(\T^2)\geq\int_{S^1}ds \frac{ \left(\int_{\ell_s} f dt\right) ^2}
{\ell_0 } = \frac{1}{\ell_0}\int_{S^1} ds \left( \length \phi(\ell_s)
\right)^2.
\]
Hence there is an $s_0$ such that $\area(\T^2)\geq \frac {\sigma}
{\ell_0} \length \phi(\ell_{s_0})^2$, so that $\length\phi(\ell_{s_0})
\leq \ell_0$.  This reduces the proof to the flat case.  Given a
lattice in $\C$, we choose a shortest lattice vector $\lambda_1$, as
well as a shortest one $\lambda_2$ not proportional to $\lambda_1$.
The inequality is now obvious from the geometry of the standard
fundamental domain for the action of $PSL(2,\Z)$ in the upper half
plane of $\C$.  \qed

We record here a slight generalization of Pu's theorem from \cite{Pu}.
The generalization follows from Gromov's inequality \eqref{tqi}.
Namely, every surface $(X,\g)$ which is not a 2-sphere satisfies
\begin{equation}
\label{pu3}
\pisys_1(\g)^2\leq \frac{\pi}{2}\area(\g),
\end{equation}
where the boundary case of equality in \eqref{pu3} is attained
precisely when, on the one hand, the surface $X$ is a real projective
plane, and on the other, the metric $\g$ is of constant Gaussian
curvature.

In both Loewner's and Pu's proofs, the area of two pointwise conformal
metrics were compared to the minimum lengths of paths in a fixed
family of curves (the family of all non contractible closed curves).
This is sometimes called the conformal-length method and analyzed in
\cite[section 5.5]{Gr1}.  C. Bavard \cite{Bav92} was able to show that
in a given conformal class, on any $n$-dimensional manifold, there is
at most one metric with maximum systolic ratio, and was able to give a
characterization of such metrics.

We conjecture a generalisation, \ref{puphi}, of Pu's inequality
\eqref{pu3}.  Let $S$ be a nonorientable surface and let
$\phi:\pi_1(S)\to \Z_2$ be a homomorphism from its fundamental group
to $\Z_2$, corresponding to a map $\hat\phi:S\to \R P^2$ of absolute degree
one.  We define the ``1-systole relative to $\phi$'', denoted
$\phisys_1(\g)$, of a metric $\g$ on $S$, by minimizing length over
loops $\gamma$ which are not in the kernel of $\phi$, \ie loops whose
image under $\phi$ is not contractible in the projective plane:
\begin{equation}
\label{cps}
\phisys_1(\g)=\min_{\phi([\gamma])\not= 0\in \Z_2}\length(\gamma).
\end{equation}
The question is whether this systole of $(S,\g)$ satisfies the
following sharp inequality, related to Gromov's inequality
$(*)_{inter}$ from \cite[3.C.1]{Gr2}, see also
\cite[Theorem~4.41]{Gr3}.
\begin{conjecture}
\label{puphi}
For any nonorientable surface $S$ and map $\hat\phi:S\to \R P^2$ of
absolute degree one, we have $\phisys_1(\g)^2\leq
\tfrac{\pi}{2}\area(\g)$.
\end{conjecture}
In the above, the ordinary 1-systole of Pu's inequality \eqref{pu3} is
replaced by the one relative to $\phi$.  The example of the connected
sum of a standard $\R P^2$ with a little 2-torus shows that such an
inequality would be optimal in every topological type $S$.
Conjecture~\ref{puphi} is closely related to the the filling area of
the circle (Conjecture~\ref{fillsphere}).  Given a filling $X^2$, one
identifies antipodal points on the boundary circle to get a
nonorientable surface $S$.  One can define a degree one map $f$ from
$S$ to $\R P^2$ by taking all points outside of a tubular neighborhood
of the original boundary circle to a point.  This gives a natural
choice of $\phi$ above which relates the two conjectures, \cf Remarks
$(e), (e')$ following \cite[Theorem 5.5.B$'$]{Gr1}.

\subsubsection{A surface of genus 2 with octahedral Weierstrass set}
\label{genus2}

The optimal systolic ratio in genus 2 is unknown.  Here we discuss a
lower bound for the optimal systolic ratio in genus 2.  The example of
M. Berger (see \cite [Example~5.6.B$'$]{Gr1}, or \cite{Be83}) in genus
2 is a singular flat metric with conical singularities.  Its systolic
ratio is 0.6666, which is not as good as the two examples we will now
discuss.

C. Bavard \cite{Bav2} and P.~Schmutz \cite[Theorem~5.2]{Sch}
identified the hyperbolic genus 2 surface with the optimal systolic
ratio among all hyperbolic genus 2 surfaces.  The surface in question
is a triangle surface (2,3,8).  It admits a regular hyperbolic octagon
as a fundamental domain, and has 12 systolic loops of length $2x$,
where $x=\cosh^{-1} (1+\sqrt{2}).$ It has $\pisys_1= 2\log\left(
1+\sqrt{2}+ \sqrt{2+2\sqrt{2}} \right)$, area ${4\pi}$, and systolic
ratio 0.7437.

This ratio can be improved by a singular flat metric $\g_{\cal O}$,
described below.  Note that it is an Hadamard space in a generalized
sense, \ie a CAT(0) space.

Start with a triangulation of the real projective plane with 3
vertices and 4 faces, corresponding to the octahedral triangulation of
the 2-sphere.  The fact that the metric is pulled back from $\R P^2$
is significant, to the extent that if we can show that an extremal
metric satisfies this, we would go a long way toward identifying it
explicitly.

Every conformal structure on the surface of genus 2 is hyperelliptic
\cite[Proposition III.7.2]{FK}, \ie it admits a ramified double cover
over the 2-sphere with 6 ramification points, called Weierstrass
points.  We take the 6 vertices of the regular octahedral
triangulation, corresponding to the Riemann surface which has an
equation of the form 
\[
y^2 = x^5-x .
\]
This produces a triangulation of the genus 2 surface $\Sigma_2$
consisting of 16 triangles.  If each of them is flat equilateral with
side $x$, then the total area is $16\left( \frac{1}{2} x^2 \sin
\frac{\pi}{3}\right)= 4\sqrt{3}x^2$.  Meanwhile, $\pisys_1= 2x$,
corresponding to the inverse image of an edge under the double cover
$\Sigma_2\to S^2$.  The systolic ratio is $\frac {(2x)^2}
{4\sqrt{3}x^2}= \frac{1} {\sqrt{3}} = 0.5773$, which is not as good as
the Bavard-Schmutz example.

Notice that the systolic loop remains locally minimizing if the angle
of the triangle with side $x$ is decreased from $\frac{\pi}{3}$ to
$\frac{\pi}{4}$ (and the space remains CAT(0)).  Therefore we replace
the flat equilateral triangles by singular flat equilateral
``triangles'', whose singularity at the center has a total angle of
$\frac{9\pi}{4}=2\pi\left( 1+\frac{1}{8} \right)$, while the angle at
each of the vertices is $\frac{\pi}{4}$.  Its barycentric subdivision
consists of 6 copies of a flat right angle triangle, denoted $\rat$,
with side $\frac{x}{2}$ and adjacent angle $\frac{\pi}{8}$.  Notice
that this results in a smooth ramification point (with a total angle
of $2\pi$) in the ramified cover.  We thus obtain a decomposition of
$\Sigma_2$ into 96 copies of the triangle~$\rat$.  

The resulting singular flat metric $\g_{\cal O}$ on $\Sigma_2$ (here
${\cal O}$ stands for ``octahedron'') has 16 singular points, with a
total angle $\frac{9\pi}{4}$ around each of them, 4 of them over each
center of a face of the $\R P^2$ decomposition.  This calculation is
consistent with M.~Troyanov's \cite{Tro} Gauss-Bonnet formula
$\sum_{\sing}\alpha(\sing)=2\gen-2$, where $\gen$ is the genus, where
the cone angle at singularity $\sing$ is $2\pi(1+\alpha(\sing))$.

Dually, the metric $\g_{\cal O}$ can be viewed as glued from six flat
regular octagons, centered on the Weierstrass points.  The
hyperelliptic involution is the 180 degree rotation on each of them.
The 1-skeleton projects to that of the inscribed cube in the 2-sphere.
The systolic ratio of the resulting metric $\g_{\cal O}$ is
\begin{equation}
\begin{array}{rcl}
\label{37}
\displaystyle \frac{\left( \pisys_1(\g_{\cal O})\right)^2} {\area
(\g_{\cal O})} & = & \displaystyle \frac{(2x)^2}{96\;\area\left(\rat
\right)} _{\phantom{a\choose b}} \\ & = & \displaystyle \frac {x^2}
{24\left ( \frac{1}{2} \left( \frac{x} {2} \right)^2 \tan\frac
{\pi}{8} \right)}_{\phantom{a\choose b}} \\ & = & \displaystyle
\left(3\sqrt{3-2\sqrt{2}}\right)^{-1^{\phantom{a}}} \\ & = & 0.8047.
\end{array}
\end{equation}
The genus 2 case is further investigated in \cite{Ca3}.

\subsubsection{Gromov's area estimates for surfaces of higher genus}
\label{highergenus}

The earliest work on isosystolic inequalities for surfaces of genus
$\gen$ is by R. Accola \cite{Ac60} and C. Blatter \cite{Bl61}.  Their
bounds on the optimal systolic ratio went to infinity with the genus,
\cf \eqref{bs}.  J. Hebda \cite{He82} and independently Yu.~Burago and
V.~Zalgaller \cite{B-Z80, B-Z88} showed that for $\gen>1$ the optimal
systolic ratio is bounded by 2 (which is not as good as \eqref{pu3},
which came later).  M.~Gromov~\cite[p.~50]{Gr1} (\cf \cite[Theorem 4,
part (1)]{Ko}) proved a general estimate which implies that if
$\Sigma_\gen$ is a closed orientable surface of genus $\gen$ with a
Riemannian metric, then
\begin{equation}
\label{1.1}
\frac {\pisys_1(\Sigma_\gen)^2} {\area(\Sigma_\gen)} < \frac {64}
{4\sqrt{\gen} + 27} .
\end{equation}
Thus, as expected, the optimal systolic ratio goes to $0$ as the genus
goes to infinity (see \eqref{55} for the correct asymptotic behavior).

Another helpful estimate is found in \cite[Corollary~5.2.B]{Gr1}.
Namely, every aspherical closed surface $(\Sigma,\g)$ admits a metric
ball $B=B_p\left(\tfrac{1}{2}\pisys_1(\g)\right) \subset \Sigma$ of
radius $\tfrac{1}{2}\pisys_1(\g)$ which satisfies
\begin{equation}
\label{tqi}
\pisys_1(\g)^2 \leq \frac{4}{3}\area(B).
\end{equation}
Furthermore, whenever a point $x\in \Sigma$ lies on a two-sided loop
which is minimizing in its free homotopy class, the metric ball
$B_x(r)\subset \Sigma$ of radius $r\leq \frac{1}{2}\pisys_1(\g)$
satisfies the estimate $2r^2 < \area(B_x(r))$.

\subsubsection{Simply connected and noncompact surfaces}
\label{2-sphere}

One of the motivating questions for this section is from \cite
[p.~135] {Gr1}, see also problem~87 in \cite{Ya} (or \cite{S-Y}).

\begin{question}
\label{primary}
For an $n$-dimensional compact manifold $X$, is there a constant
$C(X)$ such that for every Riemannian metric $\g$ on $X$, we have
\begin{equation}
\label{que}
\vol(\g)\geq C(X) L(\g)^n,
\end{equation}
where $L(\g)$ is the length of the shortest nontrivial closed
geodesic.
\end{question}
This is still open for many manifolds $X$, \eg for $X=S^n, n\geq 3$.
One could also ask the (stronger) question whether $C(X)$ depends only
on $n$.  Since $\pisys_1(\g)\geq L(\g)$, upper bounds on the optimal
systolic ratio give upper bounds on the constant $C(X)$ in
\eqref{que}.  Thus we have a positive answer for all closed surfaces
except for the two-sphere, $S^2$, which of course can have no
 nontrivial systolic inequalities.  However, Question~\ref{primary} does have
an affirmative answer in this case.  It was shown in \cite{Cr9} that
every metric $\g$ on $S^2$ satisfies the bound
\begin{equation}
\label{214}
\frac 1 {31^2}L(\g)^2\leq area(\g) .
\end{equation}
Another estimate in that paper
was $L(\g)\leq 9 D(\g),$ where $D(\g)$ represents the diameter.
Neither of these constants are best possible, and both have been
improved recently \cite{Ma,NR,Sab1}.  The best known bounds are
$L(\g)\leq 4 D(\g)$ and $\frac 1 {64}L(\g)^2\leq area(\g)$.
The conjectured best constant $C(S^2)$ in \eqref{que} for the
2-sphere (suggested to the first author by E. Calabi) is $\frac 1
{2\sqrt{3}}$, attained by the singular metric obtained by gluing two
equilateral triangles along their edges.

A natural way to find closed geodesics on a non simply connected
closed Riemannian manifold is to look for the shortest curve in a
nontrivial homotopy class, as we have done in our consideration of
$sys\pi_1(\g)$.  However, when $\pi_1(X)$ is trivial, the standard
technique is to used minimax arguments on nontrivial families of
curves.  For example, G. Birkhoff \cite{Bi27} considered 1-parameter
families of closed curves starting and ending in point curves, which
pass over $S^2$ in the sense that the induced map from $S^2$ to $S^2$
has nonzero degree.  He found a closed geodesic on $S^2$ by taking
minimum over all these families of the maximum length curve in the
family.  These minimax geodesics are stationary in the sense that they
are critical points (but not necessarily local minima) for the length
functional on the space of curves.  Thus one can find a short closed
geodesic by finding such a family where every curve in the family is
short.  This idea played an important part in the proof of \eqref{214}
above.

Instead of homotopy classes we could consider homology classes
(leading to the notion of $sys_1(X)$ \cf \eqref{51}).  That is, we
could look at 1-cycles of minimal mass in nontrivial homology classes.
When $H_1(X)$ is trivial one can instead use a minimax method on
1-cycles, similar to the one described in the previous paragraph, to
get stationary 1-cycles (also see section \ref{fillradvol}).  In some
important cases these also turn out to be closed geodesics but in any
event they are natural objects.  The basis for this is the work
Almgren \cite{Al62} and Pitts \cite [Theorem 4.10] {Pi}, who get
stationary varifolds (in all dimensions) via minimax techniques.  The
case of 1-cycles (where things are easier) was exploited by Calabi and
Cao \cite{CaCa92} in their proof that the shortest closed geodesic on
a convex surface is simple.  They use the fact that on $S^2$ this
minimax technique produces a closed geodesic.  The estimates of
A. Nabutovsky, R.~Rotman, and S. Sabourau \cite{NR,Sab1} improving
\eqref{214} mentioned above exploited these techniques, \cf
\eqref{nrs}.

Although, as mentioned above, Question \ref{primary} is still open for
general metrics on $S^n$, if one considers only convex hypersurfaces
of $\R^{n+1}$ then such a result was shown independently in \cite{Tr}
and \cite{Cr9}.  The sharp constants are still not known in this case.

Finally one can ask Question \ref{primary} for noncompact surfaces and
complete metrics (of finite area).  In fact, the question has a
positive answer for all surfaces.  Most of the cases were dealt with
as a special case of \cite [Theorem 4.4A] {Gr1}, while the other cases
(the plane and the cylinder) were dealt with in \cite{Cr9}.

\subsubsection{Optimal surfaces, existence of optimal metrics, fried eggs}
\label{optimal}

In the case of $\T^2$ or $\R P^2$, we saw in \eqref{(1.1)} and
\eqref{pu3} that there is a particular smooth metric which achieves
the optimal systolic ratio.  In general this will not be the case
unless one admits metrics with singularities.  In fact, due to
Gromov's compactness theorem \cite[theorem 5.6.C$'$] {Gr1}, there will
always be a singular metric achieving the maximal systolic ratio.
Such metrics (as well as those that have locally maximal systolic
ratio) are called ``extremal isosystolic metrics''.  See \cite{Gr1} or
\cite{Ca96} for a description of these singular metrics.

For the Klein bottle, C. Bavard \cite{Bav86} (also see \cite{Sa})
found the maximal systolic ratio:
\begin{equation}
\label{2p}
\pisys_1(\g)^2\le\frac{\pi}{2\sqrt{2}}\area(\g).
\end{equation}
He also identified the extremal isosystolic metric achieving it.  This
metric is singular, and is built out of two m\"obius strips of
constant curvature $+1$ (where the central curve has length $\pi$, and
width is $\frac \pi 2$).  They are then glued together along their
boundaries.  The metric is singular since the geodesic curvatures of
both boundaries point outward.

An extensive study of such extremal metrics was undertaken by
E. Calabi \cite{Ca96}.  He derived a number of properties of such
extremal metrics using variation methods. For points in such a space
at least two systoles must pass through each point.  (We slightly
abuse terminology by referring to noncontractible curves that have
length equal to $sys\pi _1$ ``systoles''.) One consequence of his
analysis is that any such extremal metric must be flat at points where
exactly two systoles pass through each point in a neighborhood.
Furthermore, these systoles must intersect orthogonally.

The study of solutions to the Euler - Lagrange equations for this
problem was futher taken up by R. Bryant in \cite{Br96}.


Some very interesting examples of piecewise flat (singular) metrics in
genus 3 that satisfy these criteria were presented in \cite{Ca96}.
They do not have the same systolic ratio, one has ratio $\frac 9
{8\sqrt{3}}\sim .6495$ and the other $\frac 8 {7\sqrt{3}}\sim .6598$,
but Calabi stated in that paper that these were probably both local
maxima (or at least stationary points) of the systolic ratio.  Today
\cite{Ca03} Calabi feels that the example with the smaller systolic
ratio (a ramified 4-fold cover of the octahedron with ratio $\frac 9
{8\sqrt{3}}$) is probably not locally a maximum.  The example with the
larger systolic ratio is still conjectured to be not only a local
maximum but the best possible.  Calabi \cite{Ca03,Ca3} also has an
example on a surface of genus 4 of a piecewise flat (with conical
singularities) metric (having a symmetry group of order 120) with a
systolic ratio of $\frac {(18\sqrt 7)^2} {1440\sqrt 7} = \frac{9\sqrt
7}{40}\approx 0.5953$.  It is built out of 60 flat rhombi with angles
arccos(1/8) and arccos(-1/8).  It is suspected that it is close to the
optimum systolic ratio but it is not optimum since there are points
through which only one systole passes.

An important local solution to the Euler-Lagrange equations for
extremal isosystolic metrics, the ``Fried Egg'', was also described in
\cite{Ca96}.  This example is a (non flat) metric on a disk whose
boundary is a regular hexagon all of whose angles are $\frac \pi 2$.
It arises in consideration of the problem of finding the least area
Riemannian metric in a hexagon (having the symmetries of a hexagon -
i.e.  the dihedral group of order 12 generated by reflections) such
that each point on a given side has distance exactly 2 to the opposite
side.  (The above additional symmetry assumption did not appear in the
original paper but should be included \cite{Ca03}.)  A solution to
related problems for an octagon and a triangle might be useful in
further increasing the systolic ratio in the example of section
\ref{genus2} by replacing the flat octagons (in the dual picture) or
the singular triangles with such fried eggs.

The example (piecewise flat with conical singularities) of Calabi
above in genus 3 with systolic ratio $\frac 8 {7\sqrt{3}}$ is the only
candidate that exists for metrics on surfaces of higher genus
achieving the maximal systolic ratio.  In contrast, Sabourau in
\cite{Sab02c} shows that in genus 2, no flat metric with conical
singularities (such as the example of section \ref{genus2}) can
achieve the maximal systolic ratio.


\subsection{Gromov's Filling Riemannian Manifolds}
\label{gromov}

In this section, we discuss some of the main results in \cite{Gr1}.
The reader should look at that paper for many interesting results in
this area.

\subsubsection{Filling radius + main estimate}
\label{fillradvol}

Let $X^n$ be a closed, smooth, $n$-dimensional manifold with a metric
$d_X$ (not necessarily Riemannian).  Let $A$ be a coefficient ring
(either $\Z$ or $\Z/2\Z$).  In \cite{Gr1}, Gromov introduced the
notion of the filling radius, $\fillrad(X,A)$, with respect to $A$, of
$(X,d_X)$.  There is a natural strong isometric embedding (in the
metric space sense!) $i:X\to L^\infty X$, defined by
$i(x)(\cdot)=d_X(x,\cdot)$.  The filling radius is the infimum of $r$
such that $i(X)$ bounds in the tubular $r$-neighborhood
\[
T_r(i(X)) \subset L^\infty X,
\]
in the sense that the image $i([X])$ of the fundamental class vanishes
in $H_n(T_r(i(X));A)$.  

The only Riemannian manifolds for which the precise value of the
filling radius is known \cite{Ka0} are spheres and real projective
spaces of constant curvature $K$, as well as a single additional case
of $\C P^2$ \cite{Katz1}.  Thus, $\fillrad(S^n)= \frac{1}{2} \arccos
\left( - \frac {1} {n+1} \right)K^{-\frac{1}{2}}$.  Meanwhile,
$\fillrad (\R P^n)= \frac{1}{3}\diam(\R P^n)$, and thus ``round'' real
projective spaces are extremal for the optimal inequality
$\fillrad(X)\leq \frac{1}{3} \diam(X)$, valid for all Riemannian
manifolds \cite{Ka0}.  Partial results in the direction of calculating
the filling radius for other two-point homogeneous spaces were
obtained in \cite{Katz0,Katz1,Katz2}.  An optimal inequality for the
filling radius appears in \cite{W}.

One of the fundamental estimates that allows universal inequalities in
higher dimensions is the following theorem, due to Gromov
\cite[Section 1.2, Main Theorem]{Gr1}.

\begin{theorem}
\label{filradvol}
For any closed Riemannian $n$-manifold $X$ we have
$$\vol(X)^{\frac 1 n}> \left( (n+1)(n^n)\sqrt{(n+1)!} \right)^{-1}
\fillrad(X).$$
\end{theorem}

Many of the higher dimensional universal inequalities use this
estimate along with an estimate for the filling radius.  The first of
these was Gromov's systolic theorem (below) for essential manifolds.

\begin{definition}
\label{(4.4)}
The manifold $X^n$ is called {\em essential} (over $A$) if it admits a
map, $F:X\to K$, to an aspherical space $K$ such that the induced
homomorphism $F_*:H_n(X,A)\to H_n(K,A)$ sends the fundamental class
$[X]$ to a nonzero class: $F_*([X])\not= 0 \in H_n(K,A)$.
\end{definition}
Here we must take $A=\Z/2\Z$ if $X$ is nonorientable.  In
\cite[Theorem 0.1.A]{Gr1} Gromov showed the following.

\begin{theorem}
\label{2.5}
Assume $X$ is essential over $A$.  Then every Riemannian metric $\g$
on $X$ satisfies the inequality
$$\fillrad(\g)\geq \frac 1 6 sys\pi_1(\g).$$
\end{theorem}
Combined with Theorem \ref{fillradvol}, this yields the inequality
\begin{equation}
\label{sys1}
\pisys_1(\g)^n < \left(6(n+1)(n^n)\sqrt{(n+1)!}\right)^n \vol_n(\g).
\end{equation}
This theorem generalizes both the Loewner \eqref{(1.1)} and the Pu
\eqref{pu3} theorems to higher dimensions (with non-sharp constants)
since both $\T^n$ and $\R P^n$ are essential, and in particular
answers Question \ref{primary} in the affirmative for these spaces.
Since the proof of the theorem is the model for most other estimates
of $\fillrad(X)$ we give it now.

\proof We will represent a filling as continuous map $\sigma:\Sigma
\to T_r(i(X))$ from an $(n+1)$-dimensional simplicial complex $\Sigma$
such that the restriction $\sigma|_{\partial \Sigma}:\partial
\Sigma\to i(X)$ represents a generator in $H_n(X;A)$.  We note that
for any fixed $\epsilon>0$, by taking barycentric subdivisions as
needed, we may assume that the $\sigma$-image of any simplex has
diameter less than $\epsilon$ in $L^\infty(X)$.  We note that there
can be no continuous map $f:\Sigma\to K$ which agrees with
$F\circ\sigma$ on $\partial\Sigma$, since $F\circ \sigma|_{\partial
\Sigma}$ represents $F_*[X]$ (so is not a boundary in $K$ by the
hypothesis of Definition \ref{(4.4)}).  Since $F$ is defined on $X$,
we can also think of it as defined on $i(X)\subset L^\infty (X)$
because $i$ is an embedding.
 
The proof will be by contradiction.  We assume that $\fillrad(X)<\frac
1 6 sys\pi_1(X)$, let $0<\epsilon$ be such that
$2\fillrad(X)+3\epsilon<\frac 1 3 sys\pi_1(X)$.  Take a filling
$\sigma:\Sigma\to T_r(i(X))$ (for $r=\fillrad(X)+\epsilon$) where the
image of the simplices have diameter less than $\epsilon$, and show
that $F\circ\sigma:\partial\Sigma\to K$ extends to a continuous
$f:\Sigma\to K$ giving the desired contradiction.

We construct the extension by first choosing, for every vertex $v$ of
$\Sigma-\partial \Sigma$, a point $x(v)\in X$ such that
$d(i(x(v)),\sigma(v))<r$ (which we can do since $\sigma(\Sigma)\subset
T_r(i(X))$).  (For $v\in \partial \Sigma$ just take $x(v)=\sigma(v)$)
Now we define $f(v)\equiv F(x(v))$.  Let $e$ be an edge of a simplex
in $\Sigma$ with endpoints $v_1$ and $v_2$.  The strongly isometric
nature of the imbedding $X\to L^\infty X$ and the triangle inequality
in $L^\infty X$ imply that
\[
d(x(v_1),x(v_2))=d(i(x(v_1)),i(x(v_2)))\leq 2r+\epsilon<\frac 1 3
sys\pi_1(X).
\]
Choose a shortest path $x(e)(t)$ from $x(v_1)$ to $x(v_2)$ in $X$ and
let $f(e(t))=F(x(e)(t))$.  This extends $f$ to the one-skeleton of
$\Sigma$.  Now for each two-simplex of $\Sigma$ with edges $e_1$,
$e_2$, and $e_3$, the closed curve formed by the $x(e_i)$ has length
less than $sys\pi_1(X)$.  Therefore it can be contracted to a point in
$X$.  Hence we can define a map of the two-simplex to $X$, and then
composing with $F$ to $K$ which agrees with $f$ on the edges.  Thus we
can extend $f$ to the two-skeleton.  Now since $K$ is aspherical,
there is no obstruction to extending $f$ to the rest of $\Sigma$.
This gives the desired contradiction.  \qed

In general, proofs of lower bounds on the filling radius take this
form.  That is, one assumes the filling radius is small, takes a nice
filling, and uses it to construct some map (usually skeleton by
skeleton) that one knows does not exist.

The following theorem from [Gro2, 3.$C_1$] can be thought of as a
generalisation of Loewner's inequality \eqref{(1.1)}, and a
homological analogue of \eqref{sys1}, see also Gromov's sharp stable
inequality \eqref{(1.4)}.  Let $X$ be a smooth compact $n$-dimensional
manifold.  Assume that there is a field $F$ and classes $\alpha_1,
\ldots, \alpha_n \in H^1(X,F)$ with a nonvanishing cup product
$\alpha_1 \cup \ldots \cup \alpha_n \not=0$.  Then every Riemannian
metric $\g$ on $X$ satisfies the inequality
\begin{equation}
\label{525}
\sys_1(\g)^n \le C_n \vol_n(\g),
\end{equation}
where $C_n$ is a constant depending only on the dimension.

In section \ref{2-sphere} we introduced the notion of a stationary
1-cycle.  These are usually found by the minimax techniques of Almgren
and Pitts (see \cite{Al62} and \cite [Theorem 4.10] {Pi} ).  A
stationary one cycle is built out of finitely many geodesic segments
and the mass is just the sum of the lengths of the segments.  For a
given Riemannian metric $\g$ we will let $m_1(\g)$ be the minimal mass
of a stationary 1-cycle.  Since a closed geodesic is a stationary
1-cycle we have $L(\g)\geq m_1(\g)$.

Recently these techniques have been combined with the Filling
techniques to get isosystolic type estimates in all dimensions for all
manifolds.  In \cite{N-R2} Nabutovsky and Rotman showed that
\begin{equation}
\label{nrs}
m_1(\g)\leq 2(n+2)!\fillrad(\g)
\end{equation}
and hence, by Theorem \ref{filradvol}, we have $m_1(\g)^n\leq
C(n)\vol(\g)$ (for an explicit constant $C(n)$).

In \cite{Sab1} Sabourau gives a lower bound on the Filling radius of a
generic metric on the two-sphere in terms of the minimal mass of a
1-cycle of index 1.  The minimal mass is achieved by either a simple
closed geodesic or a figure 8 geodesic.  The advantage of this
estimate is that short geodesics around thin necks (which, in nearby
generic metrics, will have index 0 as 1-cycles) can be ignored to give
better bounds on the area.

\subsubsection{Filling volume and chordal metrics}
\label{42}
In \cite{Gr1} Gromov also introduced the notion of Filling Volume,
$FillVol(N^n,d)$, for a compact manifold $N$ with a metric $d$ (here
$d$ is a distance function which is not necessarily Riemannian).  For
the actual definition one should see \cite{Gr1}, but it is shown in
\cite{Gr1} that when $n\geq 2$
\begin{equation}
\label{fillvoldef}
FillVol(N^n,d)=\inf_\g{\vol(X^{n+1},\g)}
\end{equation}
where $X$ is any fixed manifold such that $\partial X=N$ (one can even
take $X=N\times [0,\infty)$), the infimum is taken over all Riemannian
metrics $\g$ on $X$ for which the boundary distance function is $\geq
d$.  In the case $n=1$, the topology of the filling $X^2$ could affect
the infimum, as is shown by example in
\cite[Counterexamples~2.2.B]{Gr1}.

In fact, the filling volume is not known for any Riemannian metric.
However, Gromov does conjecture the following in \cite{Gr1},
immediately after Proposition 2.2.A.
\begin{conjecture}
\label{fillsphere}
$Fillvol(S^n,can)=\frac 1 2 \omega_{n+1}$, where $\omega_{n+1}$
represents the volume of the unit $(n+1)$-sphere.
\end{conjecture} 
This is still open in all dimensions.  In the case that $n=1$, as we
pointed out earlier, this is closely related to
Conjecture~\ref{puphi}.  The filling area of the circle of length
$2\pi$ with respect to the simply connected filling (by a disk) is
indeed $\tfrac{\pi}{2}$, by Pu's inequality \eqref{pu3} applied to the
projective plane obtained by identifying opposite points of the
circle.

In many cases it is more natural to consider ``chordal metrics'' than
Riemannian metrics when discussing Filling volume.  Consider a compact
manifold $X$ with boundary $\partial X$ with a Riemannian metric $\g$.
Then there is a (typically not Riemannian) metric $d_\g$ on $\partial
X$, where $d_\g(x,y)$ represents the distance in $X$ with respect to
the the metric $\g$, \ie the length of the $\g$-shortest path in $X$
between boundary points.  We call $d_\g$ the {\em chordal} metric on
$\partial X$ induced by $\g$.

Sharp filling volume estimates for chordal metrics are related to a
universal length versus volume question for a given compact manifold
$X$ with boundary $\partial X$.  This question compares the volumes,
$\vol(\g_0)$ and $\vol(\g_1)$, of two Riemannian metrics $\g_0$ and
$\g_1$ on $X$ if we know that for every pair of points $x,y \in
\partial X$ we have $d_{\g_0}(x,y)\leq d_{\g_1}(x,y)$.  Of course,
without some further assumptions on the metrics (such as some
minimizing property of geodesics) there is no general comparison
between the volumes.  For a fixed $\g_0$ (and $n\geq 2$) the statement
that for all such $\g_1$ we have $\vol(\g_0)\leq \vol(\g_1)$ is just
the statement that $FillVol(\partial X,d_{\g_0})=\vol(\g_0)$.  Note
that the standard metric on the $n$-sphere is just the chordal metric
of the $(n+1)$-dimensional hemisphere that it bounds.  So Conjecture
\ref{fillsphere} is in fact such a question.

In the computation of the Filling volume via formula
\eqref{fillvoldef} above when $d$ is the chordal distance function of
some Riemannian manifold $(X^{n+1},\partial X,\g_0)$ ($n\geq 2$) one
can not only fix the topology of $X^{n+1}$ but also restrict to
metrics $\g$ such that the Riemannian metrics on $\partial X$ gotten
by restricting $\g$ and $\g_0$ to $\partial X$ are the same (see
\cite{Cr1}).

The filling volume is known for some chordal metrics.  Gromov in
\cite{Gr1} proved this for $X^{n+1}$ a compact subdomain of $\R^{n+1}$
(or in fact for some more general flat manifolds with boundary). The
minimal entropy theorem of Besson, Courtois, and Gallot
\cite{B-C-G,B-C-G1} can be used to prove the result for compact
subdomains of symmetric spaces of negative curvature (see
\cite{Cr1}). For general convex simply connected manifolds
$(X,\partial X,\g_0)$ of negative curvature, there is a $C^3$
neighborhood in the space of metrics such that any metric $\g_1$ in
that neighborhood with $\g_1|_{\partial X}=\g_0|_{\partial X}$ and
$d_{\g_1}(x,y)\geq d_{\g_0}(x,y)$ has $\vol(\g_1)\geq \vol(\g_0)$ and
equality of the volumes implies $\g_1$ is isometric to $\g_0$ (see
\cite{C-D-S}). This leads to the conjecture:

\begin{conjecture} 
\label{fillvol}
For any compact subdomain of a simply connected space of negative
(nonpositive?) curvature of dimension $\geq 3$, the filling volume of
the boundary with the chordal metric is just the volume of the domain.
Furthermore, the domain is the unique (up to isometry) volume
minimizing filling.
\end{conjecture}

The uniqueness part of this question has applications to the boundary
rigidity problem. In that problem one considers the case when
$d_{\g_0}(x,y)= d_{\g_1}(x,y)$ for all boundary points $x$ and $y$ and
asks if $\g_0$ must be isometric to $\g_1$. In some natural cases (see
\SGM{} below) the volumes can be shown to be equal.  For example,
Gromov's result for subdomains of Euclidean space showed that they
were boundary rigid.  Again this cannot hold in general. A survey of
what is known about the boundary rigidity problem can be found in
\cite{Cr1}. There are a few natural choices for assumptions in this
case. The most general natural assumption of this sort is \SGM. The
\SGM{} condition (which is given in terms of the distance function
$d_\g:\partial X\times \partial X\to \R$ alone) would take some space
to define precisely (see \cite{Cr2}), but loosely speaking (\ie the
definitions coincide except in a few cases) it means the following:

\begin{definition}
\label{38}
A metric is ``loosely" \SGM{} if all nongrazing geodesic segments are
strongly minimizing.
\end{definition}
By a nongrazing geodesic segment we mean a segment of a geodesic which
lies in the interior of $X$ except possibly for the endpoints.  A
segment is said to minimize if its length is the distance between the
endpoints, and to strongly minimize if it is the unique such
path. (This loose definition seems to rely on more than $d_\g$ but the
relationship is worked out in \cite{Cr2}.)  Examples of such
$(X,\partial X,\g)$ are given by compact subdomains of an open ball,
$B$, in a Riemannian manifold where all geodesics segments in $B$
minimize.  The only reason not to use the ``loose" definition above is
that using a definition (such as \SGM) given in terms only of $d_\g$
guarantees that if $d_0(x,y)= d_1(x,y)$ and $\g_0$ is \SGM, then
$\g_1$ will be as well. In fact, the questions and results stated here
for the \SGM{} case also hold for manifolds satisfying the loose
definition, so the reader can treat that as a definition of \SGM{} for
the purpose of this paper.

The most general result one could hope for would be of the form: 
\begin{question}
\label{bdrig}
If $\g_0$ is an \SGM{} metric on $(X,\partial X)$ and $\g_1$ is
another metric with $d_{\g_0}(x,y)\leq d_{\g_1}(x,y)$ for $x,y\in
\partial X$ then $\vol(\g_0)\leq \vol(\g_1)$ with equality of volumes
implying $\g_0$ is isometric to $\g_1$.
\end{question}
This is still very much an open question, which as stated includes the
boundary rigidity problem.

The case $\g_1=f^2(x)\g_0$ (\ie $\g_1$ is pointwise conformal to $\g_0$)
was answered positively in \cite{C-D} (also see \cite{Cr2}).

In two dimensions more is known. Recently Ivanov \cite{Iv1} considered
the case of compact metrics $\g_0$ and $\g_1$ on a disk. He assumes
that $\g_0$ is a convex metric in the sense that every pair of
interior points can be joined by a unique geodesic, and proves that if
$d_{\g_1}\geq d_{\g_0}$ then $\vol(\g_1)\geq \vol(\g_0)$. He also says
that equality in the area would imply that $d_{\g_1}= d_{\g_0}$. Of
course, going from here to isometry of $\g_0$ and $\g_1$ is the
boundary rigidity problem. This boundary rigidity problem (in 2 dimensions) is solved
\cite {Cr2,Ot2} in the case that the metric $\g_0$ has negative
curvature.

Using different methods \cite{C-D} proves a similar though somewhat
different result where the metrics $\g_0$ and $\g_1$ are both assumed
to be $SGM$ but the surfaces are not assumed to be simply connected.
The estimate comes from a formula for the difference between the areas
of two $SGM$ surfaces.  Since the surfaces are not simply connected
one needs to worry about the homotopy class of a path between two
boundary points. We will let $L_{\g}(x,y,[\alpha])$ be the length of
the $\g$-shortest curve from $x$ to $y$ in a homotopy class $[\alpha]$
of curves from $x$ to $y$.  ${\cal A}$ will represent the space of such 
triples $(x,y,[\alpha])$ The formula is:
\begin{equation}
\label{}
area(\g_1)-area(\g_0)=\frac 1 {2\pi}\int_{\cal A} \left( L_{\g_1}
^{\phantom{a\over b}} (x,y, [\alpha])- L_{\g_0}(x,y,[\alpha]) \right)
\big(\mu_{\g_1}+\mu_{\g_0}),
\end{equation}
where the measures $\mu_{\g_i}$ on ${\cal A}$ are the push forwards of the
standard Liouville measure on the space of geodesics.  Since for each
pair $(x,y)$ there is at most one geodesic segment in each metric from
$x$ to $y$, we see that for most $(x,y,[\alpha])$ there will be
neither a $\g_1$-geodesic nor a $\g_2$-geodesic segment from $x$ to
$y$ in the class $[\alpha]$ and hence the term for $(x,y,[\alpha])$
will contribute nothing to the integral (see \cite{C-D} for
details). This formula easily leads to a result relating lengths of
paths to area.

Another powerful result with filling volume consequences is the
Besicovitch lemma which was exploited and generalized in \cite[section
7]{Gr1}.  It says that for any Riemannian metric on a cube, the volume
is bounded below by the product of the distances between opposite
faces.  Equality only holds for the Euclidean cube.

\subsubsection{Marked length spectrum and volume}
\label{mls}

The question for compact manifolds $N$ without boundary corresponding
to Question~\ref{bdrig} involves the marked length spectrum. The
marked length spectrum for a Riemannian metric $\g$ on $N$ is a
function, $\MLS_\g:{\cal C}\to \R^+$, from the set $\cal C$ of free
homotopy classes of the fundamental group $\pi_1(N)$ to the
nonnegative reals. For each $\langle \gamma\rangle\in {\cal C}$,
$\MLS_\g(\langle \gamma\rangle)$ is the length of the shortest curve
in $\langle \gamma\rangle$ (always a geodesic). We consider two
Riemannian metrics $\g_0$ and $\g_1$ on $N$ such that
$\MLS_{\g_1}(\langle \gamma\rangle)\geq \MLS_{\g_0}(\langle
\gamma\rangle)$ for all free homotopy classes $\langle \gamma\rangle$
(we then say $\MLS_{\g_1}\geq \MLS_{\g_0}$) and ask if $\vol(\g_1)$
must be greater than or equal to $\vol(\g_0)$. Again this is hopeless
without further assumptions.

The natural setting for this is in negative curvature. Here there are
lots of closed geodesics but exactly one for each free homotopy class
(achieving the minimum length in that class). The following was
conjectured in \cite{C-D-S}:

\begin{conjecture}
For two negatively curved metrics, $\g_0$ and $\g_1$, on a manifold
$N$ the inequality $\MLS_{\g_1}\geq \MLS_{\g_0}$ implies $\vol(\g_1)\geq
\vol(\g_0)$. Furthermore, $\vol(\g_1) = \vol(\g_0)$ if and only if
$\g_0$ and $\g_1$ are isometric.
\end{conjecture}

This conjecture was proved in dimension 2 in \cite{C-D}.  The higher
dimensional version of the above was shown when $\g_0$ and $\g_1$ are
pointwise conformal.  This used ideas developed in \cite{Bo1} and \cite{Si}.
 
A consequence of this conjecture is the conjecture for negatively
curved $\g_0$ and $\g_1$ that $\MLS_{\g_1}= \MLS_{\g_0}$ would imply
$\g_0$ is isometric to $\g_1$.  See \cite{Cr1} for a survey of this
problem and \cite{Ba,C-F-F,Cr2,Ot1} for results (stronger than the conjecture) 
in two dimensions. 
However, it is not even known if $\MLS_{\g_1}= \MLS_{\g_0}$
implies $\vol(\g_1)=\vol(\g_0)$.  In fact this is itself an important
open question.  Hamenst\"adt \cite{Ha2} (also see \cite{Ha1}) proved
this in the special case that $\g_0$ is further assumed to be locally
symmetric.  This was the case that had the most important immediate
applications. For example, it allows one to drop the volume assumption
in the conjugacy rigidity result in \cite{B-C-G1} and hence to see
that if $\g_0$ is locally symmetric and $\g_1$ is a metric of negative
curvature with $\MLS_{\g_1}= \MLS_{\g_0}$, then $\g_0$ is isometric to
$\g_1$.

\subsection{Systolic freedom for unstable systoles}
\label{katz}

In \cite[p.~5]{Gr1}, M. Gromov, following M. Berger, asks the
following basic question.  What is the best constant $C$, possibly
depending on the topological type of the manifold, for which the
$k$-systolic inequality \eqref{42bis} holds?  It was shown by \MK in
collaboration with A. Suciu \cite{KS1, KS2} that such a constant
typically does not exist whenever $k>1$, in sharp contrast with the
inequalities of Loewner, Pu, and Gromov discussed in sections
\ref{2dim} and \ref{gromov}.

Similar non-existence results were obtained for a pair of
complementary dimensions \eqref{ib}, culminating in the work of
I. Babenko \cite{Bab}, \cf section \ref{43}.  These results are placed
in the context of other systolic results in the table of
Figure~\ref{fig1}.

\subsubsection{Table of systolic results, definitions}
\label{sysdefs}

Let $(X,\g)$ be a Riemannian manifold, and let $k\in \N$.  Given a
homology class $\alpha\in H_1(X;A)$, denote by $\len(\alpha)$ the
least length (in the metric $\g$) of a $1$-cycle with coefficients in
$A$ representing $\alpha$.  The {\em homology 1-systole}
$\sys_1(\g,A)$ is defined by setting
\begin{equation}
\label{51}
\sys_1(\g,A)=\min_{\alpha\not= 0\in H_1(X,A)}\len(\alpha),
\end{equation}
and we let $\sys_1(\g)= \sys_1(\g,\Z)$.  The other important case is
$\sys_1(\g,\Z_2)$.  In other words, $\sys_1(\g)$ is the length of a
shortest loop which is not nullhomologous in a compact Riemannian
manifold $(X,\g)$ with a nontrivial group $H_1(X,\Z)$.

\stm

More generally, let $\sys_k(\g)$ be the infimum of $k$-volumes of
integer (Lipschitz) $k$-cycles which are not boundaries in $X$, \cf
formula~\eqref{51}.  Note that if $X$ is orientable of dimension $n$,
then the total volume is a systolic invariant:
$\vol_n(\g)=\sys_n(\g)$.

The term ``systolic freedom'' refers to the absence of a systolic
inequality, \eg~violation of the inequality relating a single systole
to the total volume, namely 
\begin{equation}
\label{42bis}
sys_k(\g)^{\frac{n}{k}} < C \vol_n(\g),
\end{equation}
or similarly the inequality involving a pair of complementary
dimensions, $k$ and $(n-k)$, namely $sys_k(\g) sys_{n-k}(\g) < C
\vol_n(\g)$, by a suitable sequence of metrics, \cf
formula~\eqref{44}.  Thus, we say that an $n$-manifold is
$k$-systolically free if
\begin{equation}
\label{41}
\inf_\g \frac{\vol_n(\g)}{\sys_k(\g)^{n/k}} =0,
\end{equation}
and $(k,n-k)$-systolically free if
\begin{equation}
\label{ib}
\inf_\g \frac{\vol_n(\g)}{\sys_k(\g)\sys_{n-k}(\g)} =0,
\end{equation}
where the infimum is over all smooth metrics $\g$ on the manifold.
See Section \ref{43} for results in this direction.  Note that we use
the reciprocals of our convention \eqref{22} for the systolic ratio,
when discussing systolic freedom.

\subsubsection{Gromov's homogeneous (1,3)-freedom}

M. Gromov first described a $(1,3)$-systolically free family of
metrics on $S^1\times S^3$ in 1993, \cf \cite[section~4.A.3]{Gr2},
\cite[p.~268]{Gr3}.  His construction of (1,3)-systolic freedom on the
manifold $S^1\times S^3$ exhibits a family $\g_\epsilon$ of
homogeneous metrics with the following asymptotic behavior:
\begin{equation}
\label{44}
\frac{\vol_4(\g_\epsilon)}{\sys_1(\g_\epsilon) \sys_3(\g_\epsilon)}
\rightarrow 0 \hbox{ as } \epsilon\to 0.
\end{equation}

Due to the exceptional importance of this example, we present a proof.
Let $U(2)$ be the unitary group acting on the unit sphere $S^3 \subset
\C^2$.  In particular, we have the scalar matrices $e^{i\theta}I_2 \in
U(2)$.  Consider the direct sum metric on $\R\times S^3$, and let
$\alpha$ be the pullback to $\R\times S^3$ of the volume form of the
second factor.  Let ${\it Trans}(\R)$ be the group of translations of
the real line.  Let $\epsilon>0$ be a real parameter.  We will define
an infinite cyclic subgroup $\Gamma_\epsilon \subset {\it
Trans}(\R)\times U(2)$ of the group of isometries of $\R\times S^3$,
in such a way that the quotient
\[
(\R\times S^3)/\Gamma_\epsilon\simeq (S^1\times S^3, \g_\epsilon)
\] 
produces a $(1,3)$-systolically free family of metrics $\g_\epsilon$
on $S^1\times S^3$ as $\epsilon\to 0$.  Namely, we may take the
generator $\gamma_\epsilon\in\Gamma_\epsilon$ to be the element
$\gamma_\epsilon=\left(\epsilon, e^{i\sqrt{\epsilon}}I_2 \right)\in
{\it Trans}(\R)\times U(2)$.  Then the quantity $\sys_1(\g_\epsilon)$
is dominated by the second factor $e^{i\sqrt{\epsilon}}$, and is
asymptotic to $\sqrt{\epsilon}$.  Meanwhile, $sys_3(\g_\epsilon)=
\omega_3$, where $\omega_3$ the volume of the unit 3-sphere, by a
calibration argument.  Namely, every nontrivial 3-cycle $C$ in
$S^1\times S^3$ satisfies $\vol_3(C)\geq \int_C \alpha = n
\omega_3\geq \omega_3$, where $n$ is the order of $[C]\in
H_3(S^1\times S^3)=\Z$.

Note that we have a fundamental domain with closure $[0,\epsilon]
\times S^3\subset \R\times S^3$ in the universal cover.  Hence
\[
\frac{\vol_4(\g_\epsilon)}{\sys_1(\g_\epsilon)
\sys_3(\g_\epsilon)}\sim \frac{\epsilon \omega_3 }{ \sqrt{\epsilon}\;
\omega_3}= \sqrt{\epsilon} \rightarrow 0,
\] 
proving (1,3)-systolic freedom of $S^1\times S^3$.

Note that a calculation of the stable 1-systole, \cf
formula~\eqref{62}, reveals the behavior of systolic constraint
instead of freedom.  Thus, the $n$-th power of the generator
$\gamma_\epsilon\in\pi_1(S^1\times S^3)=\Gamma _\epsilon$, where
$n=\left[ \frac{2\pi} {\sqrt{\epsilon}}\right]$ is the integer part,
also contains a closed geodesic whose length is on the order of
$2\pi\sqrt{\epsilon}$, as
\[
\gamma_\epsilon^n = \left(n\epsilon, e^{i\sqrt{\epsilon} n}I_2\right)\sim
\left(2\pi \epsilon^{-\frac{1}{2}}\epsilon, e^{2\pi i}I_2\right)=
\left(2\pi \epsilon^{\frac{1}{2}}, I_2\right) .
\]
It follows that the stable norm satisfies
\[
\| \gamma_\epsilon \| \sim n^{-1}\length(\gamma_\epsilon^n) \sim
\frac{2\pi \sqrt{\epsilon}}{n} \sim \epsilon,
\]
resulting in a stable 1-systole on the order of $\epsilon$ instead of
$\sqrt{\epsilon}$, which is consistent with Hebda's inequality
\eqref{hebda}.

Remarkably, M. Freedman \cite{Fr} proves the $(1,2)$-systolic freedom
of $S^1\times S^2$ even when one works with homology with coefficients
in $\Z_2$, \ie nonorientable surfaces are allowed to compete in the
definition of the 2-systole.  He showed that
\begin{equation}
\label{4.2}
\inf_\g \frac{\vol_3(\g)}{\sys_1(\g,\Z_2)\sys_{2}(\g,\Z_2)} =0,
\end{equation}
where the infimum is over all metrics $(S^1\times S^2, \g)$, \cf
formula \eqref{51}.  His technique relies on a common starting point
with the lower bound of the Buser-Sarnak Theorem \cite{BS} (\cf
\eqref{bs}), namely the existence of Riemann surfaces of arbitrarily
high genus with a uniform lower bound for the positive eigenvalues of
the Laplacian.  Whether or not $(\R P^3,\g)$ satisfies \eqref{4.2} is
unknown.

\subsubsection{Pulling back freedom by $(n,k)$-morphisms}

The early papers in the subject \cite{Ka1, Pit} contain explicit
constructions of systolically free metrics, using geometrically
controlled surgery.  Here parts of the manifold $X^n$ look like a {\em
lower dimensional\/} manifold
\begin{equation}
\label{91}
M^{m}, \hbox{ where }\; m\leq n-1,
\end{equation}
or more concretely, like circle bundle over such an $M$, with {\it
short circle fibers}.  In particular, such constructions yield free
families of metrics on {\it products of spheres}.

Later, \cite{BabK} developed a more general framework for constructing
metrics well adapted for the systolic problem.  This was done in the
context of simplicial maps $f$, where some of the top dimensional
simplices are collapsed to walls of {\em positive codimension}.  A
positive quadratic form is obtained by ``pullback'' by $f$, suitably
interpreted.  The form can then be inflated a little to make it
definite, yielding a true smooth metric.  This can be thought of as a
generalisation of the {\it short circle fiber} construction, \cf
\eqref{91}.  The next step was to realize that metrics can be pulled
back by an even more general morphism, defined below.

Let $X$ and $Y$ be $n$-dimensional simplicial complexes.  Notions of
volume, as well as systole with coefficients in $A$, can be defined
for such objects by using piecewise smooth metrics.  

A \morphism$_{n,k}$ ``from $X$ to $Y$'' (over $A$) is a continuous map
$f:X\to W$ satisfying the following two conditions:
\begin{enumerate}
\item[(a)] ~the simplicial complex $W$ is obtained from the ``target''
$Y$ by attaching cells of dimension at most $n-1$; 
\item[(b)] the map $f$ induces a monomorphism in $k$-dimensional
homology, $H_k(X,A) \hookrightarrow H_k(W,A)$.
\end{enumerate}

Note that the positive codimension of the attached cells parallels the
positive codimension in \eqref{91} above.  

Morphisms of this type were referred to as ``topological meromorphic
maps'' in \cite{KS1}.  Such a notion encompasses on the one hand, {\it
surgeries} (\cf item \ref{2} of subsection~\ref{43}), and on the
other, {\it complex blow ups}, in suitable contexts.  The importance
of these {\it morphisms}$_{n,k}$ stems from the following observation
\cite{BabK,KS1}.

\begin{proposition}
\label{48}
If $X$ admits a \morphism$_{n,k}$ to $Y$, then the $k$-systolic
freedom of $X$ follows from that of $Y$.
\end{proposition}

In the problem of systolic freedom, the appropriate {\it objects} are
not manifolds but CW-complexes $X$, even though piecewise smooth
metrics can only be defined on simplicial complexes $X'$ homotopy
equivalent to $X$.  To establish the $k$-systolic freedom for a
general $X$, one first maps $X$, by a map inducing a monomorphism on
$H_k(X,\R)$, to a product of $b_k(X)$ copies of the loop space
$\Omega(S^{k+1})$.  The space $\Omega(S^{k+1})$ is of the rational
homotopy type of the Eilenberg Maclane space $K(\Z,k)$.  

One uses the existence of a cell structure, due to Morse, which is
{\em sparse} (in the sense that only cells in an arithmetic
progression of dimensions occur) in the homotopy type of
$\Omega(S^{k+1})$.  This reduces the problem to a single
$\Omega(S^{k+1})$, rather than a Cartesian product, and thus
eliminates the dependence on the Betti number.  The $k$-systolic
freedom of the skeleta of the latter is established by finding a
morphism$_{n,k}$ to a product of $k$-spheres, \cf Proposition~\ref{48}.

This is done by analyzing the cell structure of a suitable
rationalisation, where the cases $k=2$ and $k\geq 4$ must be treated
separately ($k$ odd is easy).  In the former case, one relies on the
{\it telescope model\/} of Sullivan, and in the latter, on the
existence of rational models with cells of dimension at most $n-k+2$
added (where $n=\dim X$) in the construction of the {\it
morphism}$_{n,k}$, \cf Proposition \ref{48} and \cite{KS1, KS2} for
details.

\subsubsection{Systolic 
freedom of complex projective plane, general freedom}
\label{43}

The 2-systolic freedom of $\C P^2$ \cite{KS1} answers in the negative
the question from \cite[p.~136]{Gr1}.  Such freedom results through
the following 4 steps.

\begin{enumerate}
\item
\label{step1}
The space $\C P^2$ admits a \morphism$_{4,2}$ to $S^2\times S^2$.
Here a 3-cell is attached along a diagonal class in $\pi_2(S^2\times
S^2)$.
\item
\label{2}
The space $S^2\times S^2$ admits a \morphism$_{4,2}$ to $S^1\times
S^1\times S^2$.  Here a pair of 2-cells is attached along the
generators of $\pi_1(S^1\times S^1 \times S^2)$, corresponding to a
surgery that allows one to pass from $S^1\times S^1 \times S^2$ to
$S^2\times S^2$.

\item
There exists a $(1,2)$-systolically free sequence of metrics $\g_j$ on
the 3-manifold $S^1\times S^2$ (\cf \cite{Ka1,Pit} and
formula~\eqref{44} above).

\item
The 2-systolic freedom of $S^1\times S^1\times S^2$ results from the
Kunneth formula, by taking Cartesian product of $\g_j$ with a circle
$S^1$ of length $\sys_2(\g_j)/ \sys_1(\g_j)$.
\end{enumerate}

Part of the procedure described above can be carried out with torsion
coefficients.  Starting with Freedman's metrics \eqref{4.2} on
$S^1\times S^2$, we can construct in this way 2-systolically free
metrics on $S^2\times S^2$ even in the sense of $\Z_2$ coefficients.
However, the 2-systolic freedom of $\C P^2$ with $\Z_2$ coefficients
is still open.  The technique described above fails over $\Z_2$
because the map in step \ref{step1} has even degree, violates
condition (b) above, and hence does not define a $\morphism_{4,2}$
over $A=\Z_2$.

The series of papers \cite{Ka1, Pit, BabK, BKS, KS1, KS2, Bab} proved
general $k$-systolic freedom, as well as freedom in a pair of
complementary dimensions, as soon as any systole other than the first
one is involved.  Thus, every closed $n$-manifold is $k$-systolically
free in the sense of formula \eqref{41} whenever $2\leq k<n$ and
$H_k(X,\Z)$ is torsionfree \cite{KS1,KS2}.  This answers the basic
question of \cite[p.~5]{Gr1}.  Similarly, every $n$-dimensional
polyhedron is $(k,n-k)$-systolically free in the sense of formula
\eqref{ib} whenever $1\leq k < n-k < n$ and $H_{n-k}(X,\Z)$ is
torsionfree \cite{Bab}.  Simultaneous $(k,n-k)$-systolic freedom for a
pair of adjacent values of $k$ is explored in \cite{Ka2}.

\subsection{Stable systolic and conformal inequalities}
\label{stablesystole}

M. Gromov showed in \cite[7.4.C]{Gr1} that to multiplicative relations
in the cohomology ring of a manifold $X$, are associated stable
systolic inequalities, satisfied by an arbitrary metric on $X$, \cf
\eqref{517}.  These inequalities follow from the analogous, stronger
inequalities for conformal systoles.  He points out that the
dependence of the constants in these inequalities on the Betti numbers
is unsatisfactory.  The invariants involved are defined in \eqref{62}
and \eqref{57}.  We first point out the following open question.

\begin{conjecture}
Every closed orientable smooth 4-manifold $X$ satisfies the inequality
\begin{equation}
\label{52}
\stsys_2(\g)^2 \leq C \vol_4(\g), \;\forall\g,
\end{equation}
for a suitable numerical constant $C$ (independent of $X$).
\end{conjecture}
The lower bound of \eqref{T} shows that one is unlikely to prove such
an inequality via conformal systoles, except in the case of the
connected sum of copies of $\C P^2$, when \eqref{52} holds with $C=6$,
\cf \cite{Ka3}.  The best one has available in general is $C=6 b_2(X)$
\cite{BanK1}.

Let $H_1(X;\Z)_\R$ denote the lattice (\ie free $\Z$-module) obtained
as the quotient of $H_1(X;\Z)$ by its torsion subgroup.  The function
$i\mapsto \len(i\alpha)$ defines a norm $\|\;\|$, called {\em stable
norm} \cite[4.10]{Fe}, \cite[4.35]{Gr3}, on $H_1(X;\Z)_\R$ by
setting
\begin{equation}
\label{stable}
\|\alpha_{_\R}\|=\lim_{i\to\infty}i^{-1}{\len(i\alpha)}.
\end{equation}
Denote by $\lambda_1(L,\|\;\|)$ the least length of a nonzero vector
of a lattice $L$ with respect to a norm $\|\;\|$.  More generally, let
$i$ be an integer satisfying $1 \le i\le rk(L)$.  The $i$-th
successive minimum $\lambda_i (L, \|\;\|)$ is the least $\lambda>0$
such that there exist $i$ linearly independent vectors in $L$ of norm
at most $\lambda$:
\begin{equation}
\label{53}
\lambda_i(L,\|\;\|)= \inf_\lambda \left\{ \lambda\in \R \mid \exists
v_1,\ldots,v_i (l.i.)\; \forall k=1,\ldots,i: \|v_k\|\leq \lambda
\right\}
\end{equation} 
The {\em stable homology $k$-systole}, denoted $\stsys_k(\g)$, is the
least norm of a nonzero element in the lattice $H_k(X;\Z)_\R$ with
respect to the stable norm $\|\;\|$:
\begin{equation}
\label{62}
\stsys_k(\g)=\lambda_1\left( H_k(X;\Z)_\R, \|\;\|\right).
\end{equation}
It can be shown that if $(X,\g)$ is an orientable surface then
$\sys_1(\g,\Z)= \stsys_1(\g)$.

The conformally invariant norm $\|\;\|_{L^n}$ in $H_1(X^n,\R)$ is by
definition dual to the conformally invariant $L^n$-norm in de Rham
cohomology.  The latter norm on $H^1(X,\R)$ is the quotient norm of
the corresponding norm on closed forms.  Thus, given $\alpha\in
H^1(X,\R)$, we set
\[
\|\alpha\|_{L^n} = \displaystyle \inf_\omega \left\{ \left. \left
( \int_X |\omega_x|^n {\it dvol}(x) \right)^{\frac{1}{n}} \right|
\omega\in \alpha \right\},
\]
where $\omega$ runs over closed one-forms representing $\alpha$.  On a
surface, or more generally for a middle-dimensional class $\alpha\in
H^p(X,\R)$ where $\dim(X)=2p$, we may write
\begin{equation}
\label{54}
\|\alpha\|_{L^2}^2 = \displaystyle  \int_X \omega \wedge *\omega,
\end{equation}
where $\omega$ is the harmonic representative of $\alpha$ and $*$ is
the Hodge star operator of the metric $\g$.  The conformal 1-systole
of $(X^n,\g)$ is the quantity
\begin{equation}
\label{57}
\confsys_1(\g)=\lambda_1\left(H_1(X,\Z)_\R, \|\;\|_{L^n}\right),
\end{equation}
satisfying $\stsys_1(\g) \leq \confsys_1(\g)\vol(\g)^{\frac{1}{n}}$ on
an $n$-manifold $(X,\g)$.

\subsubsection{Asymptotic 
behavior of conformal systole as function of $\chi(X)$}
\label{51bis}

The inequality of Accola and Blatter is valid for the conformal
systole.  Their bound was improved by P. Buser and P. Sarnak
\cite{BS}, who proved that if $\Sigma_\gen$ is a closed orientable
surface of genus $\gen$, the conformal 1-systole satisfies the bounds
\begin{equation}
\label{bs}
C^{-1}\log \gen < \sup_\g \left\{ \lambda_1\left(
H^1(\Sigma_{\gen},\Z), |\;|_{L^2}^{\phantom{a}}\right) \right\} ^2 < C
\log \gen, \;\forall \gen =2,3,\ldots
\end{equation}
where $C>0$ is a numerical constant, the supremum is over all
conformal structures $\g$ on $\Sigma_{\gen}$, and $|\;|_{L^2}
^{\phantom{a}}$ is the associated $L^2$-norm, \cf formula \eqref{54}.
An explicit upper bound of $\frac{3}{\pi}\log(4\gen +3)$ in \eqref{bs}
is provided in \cite[formula (1.13)]{BS}.  Note that by Poincar\'e
duality, $\lambda_1 \left(H^1(\Sigma_{\gen},\Z), |\;|_{L^2}
^{\phantom{a}} \right) =\lambda_1 \left(H_1 ( \Sigma _{\gen}, \Z),
|\;|_{L^2}^{\phantom{a}}\right)$.

It should be kept in mind that the asymptotic behavior of the
1-systole as a function of the genus is completely different from the
conformal length.  Indeed, M.~Gromov \cite[2.C]{Gr2} reveals the
existence of a universal constant $C$ such that we have an
asymptotically {\em vanishing} upper bound
\begin{equation}
\label{55}
\frac{\sys_1(\Sigma_\gen)^2}{ \area(\Sigma_\gen)} \le C \, \frac{(\log
\gen)^2}{\gen},
\end{equation}
whenever $\Sigma_\gen$ is a closed, orientable surface of genus
$\gen\ge 2$ with a Riemannian metric.

Inequality \eqref{bs} admits the following higher-dimensional analogue
\cite{Ka2}.  Let $n\in \N$ and consider the complex projective plane
blown up at $n$ points, $\C P^2\#n \overline{\C P}^2,$ where bar
denotes reversal of orientation, while $\#$ is connected sum.  Assume
that the Surjectivity conjecture for the period map is satisfied for
such manifolds, namely every line in the positive cone of the
intersection cone occurs as the selfdual direction of a suitable
metric.  Then the conformal 2-systole satisfies the bounds
\begin{equation}
\label{T}
C^{-1} \sqrt{n} < \sup_\g \left\{ \lambda_1\left(H^2( \C P^2\#n
\overline{\C P}^2 ,\Z), |\;|^{\phantom{a}}_{L^2}\right) \right\} ^2 <
C n,
\end{equation}
as $n\to\infty,$ where $C>0$ is a numerical constant, the supremum is
taken over all smooth metrics $\g$ on $\C P^2\#n \overline{\C P}^2$,
and $|\;|_{L^2}^{\phantom{a}}$ is the norm associated with $\g$ by
formula~\eqref{54}.

It would be interesting to eliminate the dependence of inequality
\eqref{T} on the surjectivity conjecture.  Moreover, can one improve
the lower bound in \eqref{T} to linear dependence on $n$?  Here one
could try to apply an averaging argument, using Siegel's formula as in
\cite{MH}, over integral vectors satisfying $q_{n,1}(v)=-p$.  Here one
seeks a vector $v\in \R^{n,1}$ such that the integer lattice
$\Z^{n,1}\subset \R^{n,1}$ has the Conway-Thompson behavior with
respect to the positive definite form $\SR(q_{n,1},v)$.

\subsubsection{A conjectured Pu-times-Loewner inequality}
\label{amg}

We conjecture an optimal inequality for certain 4-manifolds $X$ with
first Betti number $b_1(X)=2$ which can be thought of as a product of
the inequalities of Loewner \eqref{(1.1)} and Pu \eqref{pu3}.  The
main obstacle to its proof is the absence of a generalized Pu's
inequality \eqref{puphi} (or conjecture \ref{fillsphere} for $n=1$).
For simplicity, we first state it for the manifold $\T^2\times \R
P^2$, and later discuss the relevant topological hypothesis.

\begin{conjecture}
\label{cptl}
Assume \eqref{puphi}.  Then every metric $\g$ on $\T^2\times \R P^2$
satisfies the inequality
\begin{equation}
\label{ptl}
\pisys_1(\g)^2 \stsys_1(\g)^2\leq \frac{2}{\sqrt{3}}
\frac{\pi}{2}\vol_4(\g),
\end{equation}
while a metric which satisfies the boundary case of equality must
admit a Riemannian submersion onto a Loewner-extremal torus, with
fibers which are Pu-extremal, \ie real projective planes of constant
Gaussian curvature.
\end{conjecture}

The relevant topological hypothesis is the following: $X$ should be a
4-manifold with $b_1(X)=2$, such that moreover the universal free
abelian cover $\XX$ is essential in dimension 2, \cf \cite{KKS}, so
that in particular $X$ satisfies the nonvanishing condition
$[\XX]\not=0 \in H_2 (\XX, \Z_2)$, where $[\XX]$ is the Poincar\'e
dual of the pullback by $\overline{\AJ}_X$ of the fundamental class
with compact support of the universal cover $\R^2$ of the Jacobi torus
of $X$.

We conjecture that in the boundary case of equality, both the topology
and the metrics must be right, namely: (a) the manifold $X$ must fiber
smoothly over $\T^2$ with fiber $\R P^2$; (b) the Abel-Jacobi map
$\AJ_X: X \to \T^2$ is a Riemannian submersion with fibers of constant
volume.

\subsubsection{An optimal inequality in dimension and codimension 1}

Recently, V.~Bangert and \MK \cite{BanK1} clarified the constants
involved in Gromov's inequalities mentioned at the beginning of
section \ref{stablesystole}, and showed in particular the following.
Let $k>0$, and assume $H^k(X,\R)\not= 0$.  Let $k=\sum k_j$ and assume
that the group $H^k(X,\R)$ is spanned by cup products of classes of
dimensions $k_j$.  Then
\begin{equation}
\label{517}
\prod\nolimits_j \stsys_{k_j}(\g) \le C(k) \left(\prod \nolimits_j
b_{k_j}(X)\left(1+\log b_{k_j}(X)\right) \right) \stsys_k(\g),
\;\forall\g,
\end{equation}
where the constant $C(k)$ depends only on $k$ (and not on the metric
or the Betti numbers).  On the other hand, the optimal constants in
such inequalities are generally unknown, unlike Loewner's classical
inequality \eqref{(1.1)}.  

The following sharp inequality generalizing Loewner's is proved in
\cite[pp.~259-260]{Gr3}, based on the techniques of \cite{BI}.  Assume
that the dimension, first Betti number, and real cuplength of $X$ are
all equal to $n$.  Then
\begin{equation}
\label{(1.4)}
\stsys_1(\g)^n \le (\gamma_n)^{\frac{n}{2}}\vol_n(\g), \; \forall\g
\end{equation}
where $\gamma_n$ denotes the classical Hermite constant, while
equality in \eqref{(1.4)} is attained precisely by flat tori whose
deck transformations define the densest sphere packings in dimension
$n$.  An optimal stable systolic inequality, for $n$-manifolds $X$
with first Betti number $b_1(X,\R)$ equal to one, is due to J.~Hebda
\cite{He}:
\begin{equation}
\label{hebda}
\stsys_1(\g) \sys_{n-1}(\g)\le\vol_n(\g),
\end{equation}
with equality if and only if $(X,\g)$ admits a Riemannian submersion
with connected minimal fibers onto a circle.  An optimal inequality,
involving the conformal 1-systole, is proved in \cite{BanK2}, namely
equation~\eqref{(5.2)} below.  The new inequality generalizes
simultaneously Loewner's inequality \eqref{(1.1)}, Hebda's inequality
\eqref{hebda}, the inequality \cite[Corollary 2.3]{BanK1}, as well as
certain results of G. Paternain \cite{Pa}.  We define the
Berg\'e-Martinet constant, $\gamma'_b$, by setting
\begin{equation}
\label{(2.two)}
\gamma'_b = \sup_L\left\{ \lambda_1(L) \lambda_1(L^*) \left| L \subset
\left( \R^b,|\;| \right) \right. \right\},
\end{equation}
where the supremum is over all lattices $L$ in $\R^b$, \cf
\cite{BM,CS}.  Here $L^*$ denotes the lattice dual to $L$, while
$|\;|$ is a Euclidean norm.  The supremum defining $\gamma'_b$ is
attained, and the lattices realizing it are called {\em dual-critical}
\cite{BM}.  The following is proved in \cite{BanK2}.  Let $X$ be a
compact, oriented, $n$-dimensional manifold with positive first Betti
number $b_1(X)\ge 1$.  Then every metric $\g$ on $X$ satisfies
\begin{equation}
\label{(5.2)}
\confsys_1(\g) \sys_{n-1}(\g)\leq \gamma'_{b_1(X)^{\phantom{i}}}
\vol_n(\g) ^{\frac {n-1} {n}},
\end{equation}
where equality occurs if and only if there exists a dual-critical
lattice $L$ in Euclidean space $\R^b$ (\cf~formula~\eqref{(2.two)}),
and a Riemannian submersion of $X$ onto the flat torus $\R^b/L$, such
that all fibers are connected minimal submanifolds of $X$.


\subsection{Isoembolic Inequalities}
\label{inject}

Let $inj(X)$ represent the injectivity radius of a Riemannian
manifold $X$.  All geodesics of length $\leq inj(X)$ will thus be minimizing.
This section has to do with "isoembolic inequalities" which are estimates of volume in terms of the injectivity radius.  In the closed ball, $B_p(r)$, centered at $p$ of radius $r<\frac {inj(X)} 2$ all geodesic segments will minimize and it will be \SGM \ 
as a manifold with boundary (which we denote $S(r)$).  We will let $\omega_n$ represent the volume of the standard unit sphere (of injectivity radius and diameter
$\pi$) while $\beta_n$ will be the volume of the Euclidean $n$-ball (so that $n\beta_n=\omega_{n-1}$).  One major open question is:

\begin{conjecture}
\label{volball}
\\a) For any $r\leq \frac {inj(X)} 2$, $\vol(B(r))\geq \frac {\omega_n}
2 (\frac 2 {\pi})^n r^n$ with equality holding only if the ball is
isometric to a hemisphere.  \\b) For any $r \leq inj(X)$,
$\vol(B(r))\ge \frac {\omega_n} {\pi^n} r^n$, where equality holds if
and only if $X$ is isometric to the round sphere of injectivity radius
$r$ (\ie extrinsic radius $\frac r \pi$).
\end{conjecture}

In part a) above, the hemisphere in the equality case will be that of
a sphere of (extrinsic) radius $\frac {2r} {\pi}$ (so the hemisphere
is a metric ball of intrinsic radius $r$).  In fact a) is a stronger
conjecture than b) since by the triangle inequality $B(r)$ contains
two balls of radius $\frac r 2$ with disjoint interiors.

\subsubsection{Berger's 2 and 3 dimensional estimates}

The early results on this conjecture and the conjecture itself are due
to Berger (see \cite{Be2}).  The known results take the form of
proving inequalities $\vol(B(r))\geq c(n) r^n$ for non-sharp constants
$c(n)$.  In 2 dimensions it is not hard to get such an estimate when
$r\leq \frac {inj(X)} 2$.  Simply notice that for all $t\leq r$,
$\vol(S_p(t))\geq 4t$ (\ie the length of the boundary of $B_p(t)$ is
$\geq 4t$) since ``antipodal'' points on $S(t)$ must be at least $2t$
apart in $X$ and hence also along $S(t)$.  This says that
$\vol(B_p(r))=\int_0^r \vol(S_p(t)) dt \geq 2r^2$.  The conjecture
would say $\vol(B(r))\geq \frac 8 \pi r^2$ and hence this simple
argument leads to a constant not too far from the best possible.  The
conjecture with the sharp constant is still open.  (Using ideas in
\cite{C-D} - \cf section \ref{42} - one can show in this 2-dimensional
setting that $\vol(B_p(t))\geq 4\pi -\vol(S_p(t))$ and then use this
to get the better estimate $\vol(B_p(r))\geq \frac {16} {4+ \pi} r^2$,
but since this argument is a little involved and does not yield the
sharp constant we will not pursue it.)

An interesting example to keep in mind in the ``Mercedes-Benz'' example.  Consider three unit length line segments in the plane emanating from the origin making angles of $\frac {2\pi} {3}$.  Thin tubular neighborhoods have arbitrarily small area but every boundary point is at a bounded distance from its antipodal point (the point halfway around the boundary).  Thus one needs more than just bounded distance between ``antipodal boundary points'' to get area lower bounds.

One might hope to generalize the above simple argument for the area of 2-balls by showing:

\begin{question}
\label{volsphere}
Let $\g$ be a Riemannian metric on the $n$-sphere $S^n$ such that the
$\g$-distance between antipodal points is $\geq t$.  \\a) Is
$\vol((S^n,\g))\geq \frac {\omega_n}{\pi^n} t^n$?  \\b) Is there any
constant $c(n)$ such that $\vol((S^n,\g))\geq c(n)t^n$?
\end{question}

In the above one can either think of the standard antipodal map when
$S^n$ is embedded in $\R^n$ in the standard way, or simply take the
antipodal map to be any order two fixed point free diffeomorphism.
Metric balls of radius $\frac t 2$ satisfy this when $t\leq inj(X)$.

The easy argument we gave above for the volume of 2-balls consists of
answering a) above in the case $n=1$.  This is the only case where a)
is known.  Berger in \cite{Be2} gave an answer to b) above for $n=2$
and hence showed by integrating (as in the 2 dimensional argument)
that there is a constant $c$ such that for 3 dimensional manifolds
when $r\leq \frac {inj(X)} 2$, then $\vol(B_p(r))\geq cr^3$.

Though part a) of Question \ref{volsphere} in 2 dimensions is still
open (and very interesting), there are enough reasons to believe it to
be true that it would warrant being called a conjecture.  One approach
would be to use the uniformization theorem and conformal length techniques
(but note that the antipodal map need not be conformal).

In higher dimensions there is not enough evidence one way or the other
even for question b).  However (see \cite{Cr5}) there are constants
such that for any given metric $g$ on the 3-sphere either Question
\ref{volsphere} b) holds or Question \ref{primary} holds.  This
follows from a lower bound on the filling radius (\cf section
\ref{fillradvol}) in terms of the infimum of the length of closed
curves who link their antipodal images along with an application of
Theorem \ref{filradvol}.  In the other direction, Ivanov (see
\cite{Iv2} or \cite{Iv3}) has given examples of a sequence of metrics
on $S^3$ that Gromov-Housdorff converge to the standard metric but
whose volumes go to zero.  Although this does not give a counter
example it does show that the topological properties of the antipodal
map are important in the above question.  (Even in two dimensions long
thin cigar shapes show that simply because every point is far from
some other point the area need not be large.)

\subsubsection{Higher dimensions}

The two ways that higher dimensional isoembolic type inequalities have
been proved are via Gromov's estimate of section \ref{fillradvol} and
via an estimate of Berger and Kazdan \cite{Be-Ka}.  The Berger-Kazdan
inequality first appeared in the proof of the Baschke conjecture for
spheres (see \cite{Bes} appendices D and E).  Berger then used the result to
give the sharp estimate which holds for any compact Riemannian
manifold:
\begin{equation}
\label{bergerisoemb}
\vol(X)\geq \frac {\omega_n} {\pi^n} inj(X)^n
\end{equation}
with equality holding only for round spheres.  If you know that $X$ is
not homeomorphic to a sphere then you can do better \cite{Cr6},
showing $\vol(X)\geq c(n) inj(X)^n$ where $c(n)$ is an explicit
constant larger than $\frac {\omega_n} {\pi^n}$.

The Berger-Kazdan inequality along with Santal\'o's formula
\cite[Chapter 19]{Sa1,Sa2} was used in \cite{Cr7} to give a general
sharp isoperimetric inequality for compact Riemannian manifolds, $X$,
with boundary, $\partial X$ where all geodesic minimize (\eg \SGM{}
manifolds):
$$
\frac {\vol(\partial X)} {\vol(X)^{\frac {n-1} n}}\geq C(n)
$$ 
where the constant $C(n)$ is just $\omega_{n-1} (\frac {2}
{\omega_n})^{\frac {n-1} n}$.  Equality holds if and only if $X$ is
isometric to a hemisphere.  (The inequality has a version for all
manifolds with boundary which involves a term measuring the fraction
of geodesics that minimize to the boundary.)  One gets a (non-sharp)
version of conjecture \ref{volball} a) in all dimensions by applying
this to the balls $B_p(t)$ and integrating $t$ from $0$ to $r$:
\begin{equation}
\label{volofball}
\vol(B_p(r))\geq  C(n)^n r^n.
\end{equation}
Further, applying the isoperimetric inequality again we see:
$\vol(S_p(r))\geq C(n)^n r^{n-1}.$ The failure of this estimate to be
sharp comes from the fact that the isoperimetric inequality is not
sharp for spherical caps (\ie metric balls in the round sphere) unless
they are hemispheres.  We mentioned before (Conjecture \ref{volball}
part a) that, for $r\leq \frac {inj} 2$, we do not know the optimal
lower bound for $\vol(B_p(r))$.  Except in case of dimension 2, we
also don't know the optimal lower bound for $\vol(S_p(r))$, however
the above isoperimetric inequality will imply the sharp estimate on
$\vol(S_p(r))$ if Conjecture \ref{volball} part a is proven.

Although we are unable to prove Conjecture \ref{volball} for all balls
there are some results about the ``average'' volume of balls,
$AveVol(r)$, of a compact Riemannian manifold $X$.  By this we mean
$$AveVol(r)\equiv \frac 1 {\vol(X)}\int_X \vol(B_p(r)) dp$$ where $dp$
is the Riemannian volume form.  In \cite{Cr8} it is shown for any
$r\leq inj(X)$ that $AveVol(r)\geq \frac {\omega_n} {\pi^n} r^n$ with
equality holding only for the round sphere of injectivity radius $r$.
Thus showing that part b) of Conjecture \ref{volball} is true ``on the
average'' (and giving another proof of Berger's isoembolic inequality
\ref{bergerisoemb}).  If one further knows that $X$ has no conjugate
points, then one can improve this to $AveVol(r)\geq \beta_n r^n$,
where equality holds if and only if $X$ is flat \cite{Cr4}.  In fact,
this holds for all $r$, if one interprets $\vol(B_p(r))$ to mean the
volume of a ball in the universal cover centered at a lift of $p$.

Gromov in \cite{Gr1}, showed (in an argument similar to the one
presented in section \ref{fillradvol}) that the filling radius of a
compact Riemannian manifold is always bounded from below by a constant
times the injectivity radius.  This, along with Theorem
\ref{filradvol}, gave an alternative proof to Berger's isoembolic
inequality \ref{bergerisoemb} (albeit with a non-sharp constant).
This idea was extended in \cite{G-P} to get non-sharp universal
estimates on the volumes of balls in terms of the local geometric
contractibility function.  As a special case, it gives an alternative
proof of inequality \ref{volofball}.  Sabourau also used these filling
radius ideas in \cite{Sab2} to give improved versions of both
\eqref{bergerisoemb} and \eqref{volofball} (with worse constants)
where instead of injectivity radius he is able to substitute the
length of the shortest geodesic loop. A geodesic loop is a closed
curve which is a geodesic at all but one point.  Geodesic loops are
much more abundant than closed geodesics and the shortest one will
have length, $sgl(\g)$, which satisfies
\begin{equation}
\label{65}
2inj(\g) \leq sgl(\g) \leq L(\g)
\end{equation}
Sabourau's version of \eqref{volofball}, which holds for all $r\leq
\frac 1 2 sgl(\g)$, also used the ideas of \cite{G-P}.

\subsection{Acknowledgments}
We are grateful to E. Calabi, I. Chavel, A. Nabutovsky, and
S. Sabourau for reading a draft version of the manuscript and for
their comments.


\bibliographystyle{amsalpha}

\no
Christopher B Croke:\\
Department of Mathematics\\
University of Pennsylvania\\
Philadelphia, PA 19104-6395\\
ccroke@math.upenn.edu\\

\no 
Mikhail G. Katz\\
Department of Mathematics and Statistics\\
Bar Ilan University\\
Ramat Gan  52900, Israel\\
katzmik@math.biu.ac.il\\

\end{document}